\documentclass[reqno,11pt,a4paper]{amsart}
\usepackage{color,latexsym,amssymb,amsfonts,amsmath}
\usepackage{graphicx}

\renewcommand{\b}[1]{\mbox{\boldmath $#1$}}

\newtheorem{remark}{Remark}

\numberwithin{equation}{section}
\DeclareMathOperator*{\trace}{trace}
\DeclareMathOperator*{\argmin}{argmin}

\begin{document}
\author{Goldentayer, L.}
\address{Electrical Engineering Systems,
Tel Aviv University, 69978 - Ramat Aviv, Tel Aviv, Israel}
\email{goldlev@eng.tau.ac.il}
\author{Liptser, R.}
\address{Electrical Engineering Systems,
Tel Aviv University, 69978 - Ramat Aviv, Tel Aviv, Israel}
\email{liptser@eng.tau.ac.il}
\date{16/06/2002}
\title[On-line tracking of a smooth regression function]
{On-line tracking of a smooth regression function}
\begin{abstract}
We construct an on-line estimator with equidistant design for
tracking a smooth function from Stone-Ibragimov-Khasminskii class.
This estimator has the optimal convergence rate of risk to zero in sample size.
The procedure for setting coefficients of the estimator
is controlled by a single parameter and has a simple numerical
solution. The off-line version of this estimator allows
to eliminate a boundary layer. Simulation results are given.
\end{abstract}

\thanks{\bf This work is partially supported by a fellowship from the
Yitzhak and Chaya Weinstein Research Institute for Signal
Processing at Tel Aviv University.}
\maketitle

\section{Introduction.}

In this paper, we consider a tracking problem for smooth function $f=f(t)$, $0\le t\le T$,
under observation
\begin{equation}\label{model}
X_{in} =f(t_{in})+\sigma\xi_i,
\end{equation}
for $t_{in}=\frac{i}{n}, \;\; i=0,\ldots,n$ ($n$ is large),
where $(\xi_i)$ is a sequence of i.i.d. random variables with
$E\xi_i=0,\;\; E\xi^2_i = 1$, and $\sigma^2$ is a positive constant.
Without additional assumptions on the function $f$ it is
difficult to create an estimator even for large $n$. The filtering approach,
see Bar-Shalom and Li  \cite{B_Sh}, proposes an estimator in the form of Kalman filter
corresponding
to a stochastic model for $f$, e.g. $f$ is differentiable
$k$ times, and $k$-th derivatives of $f$ is simulated by a white noise with a certain
intensity. Since $f$ is deterministic function, the non-trivial part of such
approach is a choice of filter parameters
and asymptotic analysis of estimation risk in $n\to \infty$.
On the other hand, nonparametric statistic approach to the regression
estimation of a function $f$ assumes that $f$ belong to some limited class.
We take $f$ from the class $\Sigma(\beta,L)$
(introduced by Stone, \cite{St80}, \cite{St82} and
Ibragimov and Khasminskii, \cite{IK80},  \cite{IK81})
of $k$ times continuously differentiable functions with H\"older continuous last
derivative (here $f^{(0)}=f$, $L$ and $\alpha$ are the same for any function from the
class):
$$
\Sigma(\beta,L)=\left\{f:\begin{array}{lll}
\mbox{obeys  $k$ derivatives}, f^{(0)},f^{(1)},\ldots,f^{(k)};
\\
|f^{(k)}(t_2)-f^{(k)}(t_1)|\le L|t_2-t_1|^\alpha, \ \forall \ t_1,t_2, \
 \alpha\in(0,1];
 \\
 \beta=k+\alpha
\end{array}
\right\}.
$$
It is known
 \cite{St80}, \cite{St82}, \cite{IK80},  \cite{IK81} that
there are kernel type estimators $\widehat{f^{(j)}_n}(t)$ of
$f^{(j)}(t)$, $j=0,1,\ldots,k$ such that for a wide class of loss
functions $\mathcal{L}(*)$ and $p>1$ ($C$ is positive constant):
\begin{equation}\label{riskbound}
  \sup_{f\in\Sigma(\beta,L)}E\mathcal{L}\Big(n^{\frac{\beta-j}{2\beta+1}}
  \| \widehat{f^{(j)}_n}-f^{(j)}\|_{L_p}\Big)< C, \ j=0,1,...,k
\end{equation}
and no estimator provides a better rate of convergence to zero in
$n\to\infty$ uniformly in $\Sigma(\beta,L)$. The same rate in $n$
is valid under fixed value $t$ for the estimation risk
$E\big(\widehat{f^{(j)}_n}(t)-f^{(j)}(t)\big)^2$,
$j=0,1\ldots,k$. This rate cannot be exceeded uniformly on any
nonempty open set from $(0,T)$. Parallel to kernel type
estimators (see, e.g. \cite{IK81}), \cite{R}, \cite{P}),
Khasminskii and Liptser \cite{KhL} proposed an on-line estimator
(hereafter for brevity $t_{in}$ is replaced by $t_i$ and
$\widehat{f^{(j)}_n}$ by $\widehat{f^{(j)}}$):
\begin{equation}\label{1.4tri}
\begin{aligned}
\widehat{f^{(j)}}(t_i)&=\widehat{f^{(j)}}(t_{i-1})+
\frac{1}{n}\widehat{f^{(j+1)}}(t_{i-1})+ \frac{q_j}{
n^{\frac{(2\beta-j)}{2\beta+1}}} \big(X_i-\widehat{f^{(0)}}(t_{i-1})\big)
\\
j&=0,1,\ldots,k-1
\\
\widehat{f^{(k)}}(t_i)&=\widehat{f^{(k)}}(t_{i-1})+
\frac{q_k}{n^{\frac{(2\beta-k)}{2\beta+1}}}
\big(X_i-\widehat{f^{(0)}}(t_{i-1})\big),
\end{aligned}
\end{equation}
subject to the initial conditions
$\widehat{f^{(0)}}(0),\widehat{f^{(1)}}(0),\ldots,
\widehat{f^{(k)}}(0)$. The initial conditions are chosen as arbitrary bounded constants
independent of $n$. The parameters $q_0,\ldots,q_k$ are specifically chosen.
The vector $\mathfrak{q}$ with these entries is called the filter gain.

The rigorous result given in \cite{KhL} is formulated as:
\begin{quote}
{\it
Let filter gain $\mathfrak{q}$ be chosen such that all
roots of characteristic polynomial
\begin{eqnarray}
p^k(u,\mathfrak{q})=u^{k+1}+q_0u^k+q_1u^{k-1}+\ldots+q_{k-1}u+q_k
\label{1.5}
\end{eqnarray}
are different and have negative real parts.
Let the observation model defined in (\ref{model}), $f\in
\Sigma(\beta,L)$ and $\sigma^2>0$. Then, for the estimator given in (\ref{1.4tri})
there exist positive constants $c(\mathfrak{q})$, $C(\mathfrak{q})$ (independent of $n$)
such that for any $t_i\ge c(\mathfrak{q})n^{-\frac{1}{2\beta+1}}\log n$ the
normalized in $n$ risk obeys
\begin{equation}\label{main0}
\varlimsup_{n\to\infty}
\sup_{f\in\sum(\beta,L)}\sum_{j=0}^kE\big(f^{(j)}(t_i)-
\widehat{f^{(j)}}(t_i)\big)^2
n^{\frac{2(\beta-j)}{2\beta+1}}\le C(\mathfrak{q}).
\end{equation}
The rates $n^{-\frac{2(\beta-j)}{2\beta+1}}$, $j=0,1,\ldots,k$
cannot be improved. The boundary layer
$c(\mathfrak{q})n^{-\frac1{2\beta+1}}\log n$, where \eqref{main0} might fail, is
inevitable.
}
\end{quote}

\begin{remark}
The left side boundary layer $c(\mathfrak{q})n^{-\frac1{2\beta+1}}\log n$ is due to
on-line limitations of the above tracking system. One can readily suggest an off-line
modification with the same recursion in the backward time subject to some boundary
conditions independent of observation $X_i$'s. This modification possesses
the right side boundary layer
$
[T-c(\mathfrak{q})n^{-\frac1{2\beta+1}}\log n,T]
$
and accuracy (\ref{main0}) on $[0,T-c(\mathfrak{q})n^{-\frac1{2\beta+1}}\log n]$.
So, some combination of the forward and backward time tracking algorithms allows
(\ref{main0}) accuracy on $[0,T]$. For instance, a combination of forward time tracking
on the interval $\big[\frac{T}{2},T\big]$ and backward time tracking
on $\big[0,\frac{T}{2}\big]$ can be used.
\end{remark}

In this paper, we deal with the estimator given in (\ref{1.4tri}) and restrict
ourselves by considering $f$ from the class $\Sigma(L,k+1)$, i.e.
the class of $k$-times differentiable functions $f$ with Lipschitz continuous
$f^{(k)}(t)$.

A suitable choice of filtering gain $\mathfrak{q}$ should satisfy
multiple requirements regarding the cost function
$C(\mathfrak{q})$ and parameter $c(\mathfrak{q})$, involved in the
description of boundary layer. Moreover, a correct choice of
$\mathfrak{q}$ should guarantee that the roots of characteristic
polynomial $p^k(u,\mathfrak{q})$ are different and have negative
real parts. These requirements might contradict each other. To
avoid contradictions, we use the fact that estimator
\eqref{1.4tri} has a structure of Kalman filter. We build a Kalman
filter according to Bar-Shalom and Li \cite{B_Sh}, so
that $f^{(k)}(t)$ is generated by a white noise with intensity
$\gamma$. For each $\gamma$, we choose the Kalman gain
$\mathfrak{q}(\gamma)$ and use it for minimization of
$C\big(\mathfrak{q}(\gamma)\big)$ in $\gamma$. So, the
minimization problem of the cost function is controlled by single
parameter and allows to establish a reasonable relationship
between $C(\mathfrak{q})$ and $c(\mathfrak{q})$. Moreover, this
type of minimization automatically guarantees negative real parts
of the roots for characteristic  polynomial $p^k(u,\mathfrak{q})$.
For $k\le 4$ the roots of $p^k(u,\mathfrak{q})$ are different and
numerical verification of the same fact for $k>4$ is available.

\section{The filter gain choice}
\mbox{}

\subsection{Preliminaries}

Henceforth $\beta=k+1$. For notational convenience
we describe our problem in matrix notation.
Introduce the following matrices:
$$
a=
\begin{pmatrix}
0 & 1 & 0 & 0 & \ldots & 0
\\
0 & 0 & 1 & 0 & \ldots & 0
\\
\vdots & \vdots & \vdots & \vdots & \vdots & \vdots
\\
0 & 0 & 0 & 0 & \ldots & 1
\\
0 & 0 & 0 & 0 & \ldots & 0
\end{pmatrix}_{(k+1)\times(k+1)}, \
$$

$$
A=\begin{pmatrix}
    1 & 0 &\ldots & 0
  \end{pmatrix}_{1\times(k+1)}, \quad
\mathfrak{b}=
  \begin{pmatrix}
    0 \\
    \vdots\\
    0\\
    1
  \end{pmatrix}_{(k+1)\times 1}.
$$ Notice that by Lemma 3.1 in \cite{CKL} the  roots of
$p^k(u,\mathfrak{q})$ and eigenvalues of $(a-\mathfrak{q}A)$
coincide. In accordance with this remark, while eigenvalues of
$(a-\mathfrak{q}A)$ have negative real parts we may describe the
cost function $C(\mathfrak{q})$ in terms of the bias
$\widetilde{M}(\mathfrak{q})$ and variance $P(\mathfrak{q})$ for
tracking errors (see, (\ref{prosto})) (hereafter $^*$ is the
transposition symbol):
\begin{eqnarray*}
C(\mathfrak{q})=
\trace\Big(P\big(\mathfrak{q})+\widetilde{M}\big(\mathfrak{q}\big)
\widetilde{M}^*\big(\mathfrak{q}\big)\Big)
\end{eqnarray*}
where
$
\widetilde{M}\big(\mathfrak{q}\big)=L(a-\mathfrak{q}A)^{-1}\mathfrak{b}
$
and the matrix $P\big(\mathfrak{q}\big)$ solves the Lyapunov equation
$
\big(a-\mathfrak{q}A\big)P\big(\mathfrak{q}\big)+
P\big(\mathfrak{q}\big)\big(a-\mathfrak{q}A\big)^*
+\sigma^2\mathfrak{q}\mathfrak{q}^*=0.
$
In Section \ref{sec-4}, we select the filter gain $\mathfrak{q}(\gamma)$
from one-parameter family
$$
\Gamma=\Big\{\gamma\ge \gamma_\varepsilon>0:\mathfrak{q}(\gamma)=\frac{Q(\gamma)A^*}
{\sigma^2}\Big\}
$$
where $Q(\gamma)$ is a positive definite matrix given by the algebraic Riccati equation
$
aQ(\gamma)+Q(\gamma)a^*+\gamma^2\mathfrak{b}\mathfrak{b}^*-
\frac{Q(\gamma)A^*AQ(\gamma)}{\sigma^2}=0.
$
Finally we choose
$$
\gamma^\circ=\argmin_{\gamma\in\Gamma}C\big(\mathfrak{q}(\gamma)\big)
$$
and the filter gain $\mathfrak{q}(\gamma^\circ)$. A relevant choice of
$\gamma_\varepsilon$ allows to have an acceptable value of the constant
$c\big(\mathfrak{q}(\gamma^\circ)\big)$.

\subsection{Explicit formulae}

In Section \ref{sec-5}, we show that the cost function is expressed as:
$$
C\big(\mathfrak{q}\big)=\sigma^2\Bigg(\int_0^\infty\mathfrak{q}^*e^{(a-\mathfrak{q}A)^*t}
e^{(a-\mathfrak{q}A)t}\mathfrak{q}dt+\Big(\frac{L}{\sigma}\Big)^2
\Big[\Big(\frac{1}{q_k}\Big)^2+\sum_{j=0}^{k-1}
\Big(\frac{q_j}{q_k}\Big)^2\Big]\Bigg).
$$
Furthermore, we give the explicit structure of
$\mathfrak{q}(\gamma)$ as a function of the control parameter
$(\gamma/\sigma)^{1/(k+1)}$. Namely
\begin{eqnarray*}
&&
q_0(\gamma)=U_{00}\Big(\frac{\gamma}{\sigma}\Big)^{1/k+1}
\\
&&
q_1(\gamma)=U_{01}\Big(\frac{\gamma}{\sigma}\Big)^{2/k+1}
\\
&&.....................................
\\
&&
q_k(\gamma)=U_{0k}\Big(\frac{\gamma}{\sigma}\Big),
\end{eqnarray*}
where $U_{ij}$, $i,j=0,1,\ldots,k$ are entries of the matrix $U$ being solution
of the algebraic Riccati equation
$
aU+Ua^*+\mathfrak{b}\mathfrak{b}^*-UA^*AU=0.
$

\subsection{Example 1}

Here, we consider the tracking problem for Lipschitz continuous function $f$.
Since $k=0$, we have $a=0$, $A=1$, $\mathfrak{b}=1$, $1-U^2_{00}=0$
and so, $q_0=\frac{\gamma}{\sigma}$. Therefore,
\[
C(\mathfrak{q}(\gamma))=\frac{\sigma\gamma}{2}+\frac{L^2\sigma^2}{\gamma^2}.
\]
With $\gamma_\varepsilon<(2L)^{2/3}\sigma^{1/3}$, we have
$\gamma^\circ=(2L)^{2/3}\sigma^{1/3}$ and $\mathfrak{q}(\gamma^\circ)=
\Big(\frac{2L}{\sigma}\Big)^{2/3}$.
The following estimator is constructed (here $\widehat{f}(t_i):=
\widehat{f^{(0)}}_n(t_i)$)
$$
\widehat{f}(t_i)=\widehat{f}(t_{i-1})+\Big(\frac{2L}{n\sigma}\Big)^{2/3}
(X_i-\widehat{f}(t_{i-1})).
$$

\begin{remark}
Notice that the direct  minimization of
$
C(\mathfrak{q})=P(\mathfrak{q})+\widetilde{M}^2(\mathfrak{q})
$
with  respect to $q_0$ provides the optimal $q_0=q(\gamma^\circ)$.
For $k\ge 1$, this coincidence is not guaranteed. Under the direct
minimization of the cost function $C(\mathfrak{q})$ with respect to $\mathfrak{q}$
the eigenvalues of $(a-\mathfrak{q}^*A)$ might have nonnegative real parts.
\end{remark}

\subsection{Example 2}

\begin{figure}
\begin{center}
\includegraphics[scale=0.70]{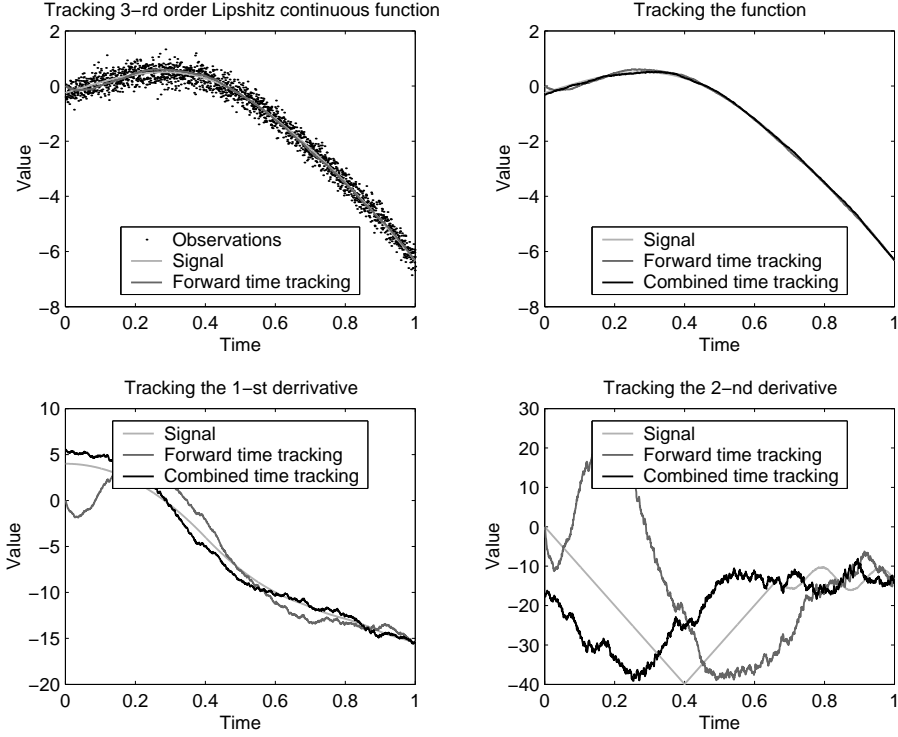}
\end{center}
\caption{Forward and backward time tracking with $n=2000$}
\label{fig2}
\end{figure}

Let us consider a numerical solution for $k=2$, i.e. $f$ is twice differentiable function. Its second
derivative is Lipschitz continuous with constant $L=100$. For $\sigma=0.25$, we find
(see Figure \ref{fig1})
$\gamma^\circ=24.533$ and
$\gamma^\circ/\sigma=98.132$. According to Table 1, $U_{00}=2$,
$U_{01}=2$, $U_{02}=1$. Hence $q_0=9.225$, $q_1=42.550$, and
$q_2=98.132$. So, the following estimator is constructed $$
\begin{aligned}
\widehat{f^{(0)}}(t_i)&=\widehat{f^{(0)}}(t_{i-1})+
\frac{1}{n}\widehat{f^{(1)}}(t_{i-1})+ \frac{9.225}{ n^{6/7}}
\big(X_i-\widehat{f^{(0)}}(t_{i-1})\big)
\\
\widehat{f^{(1)}}(t_i)&=\widehat{f^{(1)}}(t_{i-1})+
\frac{1}{n}\widehat{f^{(2)}}(t_{i-1})+ \frac{42.550}{ n^{5/7}}
\big(X_i-\widehat{f^{(0)}}(t_{i-1})\big)
\\
\widehat{f^{(2)}}(t_i)&=\widehat{f^{(2)}}(t_{i-1})+
\frac{98.132}{n^{4/7}} \big(X_i-\widehat{f^{(0)}}(t_{i-1})\big).
\end{aligned}
$$
The combination of
forward and backward tracking practically allows to eliminate the
boundary layer (see Figure \ref{fig2}).

\section{Error analysis}
\label{sec-3}
\mbox{}

In this section, we present a vector of normalized tracking errors and
derive the expression for $C(\mathfrak{q})$.

\subsection{Notation}
For notational convenience, set
$$
F(t_i)=
  \begin{pmatrix}
    f^{(0)}(t_i)
    \\
    f^{(1)}(t_i)
    \\
    \vdots
    \\
    f^{(k)}(t_i)
  \end{pmatrix},
\
\widehat{F}^n(t_i)=\begin{pmatrix}
    \widehat{f^{(0)}_n}(t_i)
    \\
    \widehat{f^{(1)}_n}(t_i)
    \\
    \vdots
    \\
    \widehat{f^{(k)}_n}(t_i)
\end{pmatrix}
$$
and introduce a diagonal matrix $C_n$ and vector $\mathfrak{q}_n$:
$$
C_n={\sf diag}\big(n^{\frac{\beta}{2\beta+1}},n^{\frac{\beta-1}{2\beta+1}},\dots,
n^{\frac{\beta-k}{2\beta+1}}\big)
\quad\text{and}\quad
\mathfrak{q_n}=\begin{pmatrix}
q_0n^{-\frac{2\beta}{2\beta+1}}
\\
q_1n^{-\frac{2\beta-1}{2\beta+1}}
\\
\vdots
\\
q_kn^{-\frac{2\beta-k}{2\beta+1}}
\end{pmatrix}.
$$
We use a vector-matrix form of estimator \eqref{1.4tri}
\begin{equation}\label{vmf}
\widehat{F}^n(t_i)=\widehat{F}^n(t_{i-1})+\frac{1}{n}a\widehat{F}^n(t_{i-1})+
\mathfrak{q}_n\big(X_i-A\widehat{F}(t_{i-1})\big)
\end{equation}
and obvious identity

\begin{equation}\label{vmf0}
F(t_i)\equiv
F(t_{i-1})+\frac{1}{n}aF(t_{i-1})+\mathfrak{b}\big(f^{(k)}(t_i)-f^{(k)}(t_{i-1})\big).
\end{equation}

\subsection{Normalized errors}

Denote
$
\delta_i=\widehat{F}(t_i)-F(t_i)
$
and introduce normalized error
$
\Delta_i=C_n\delta_i.
$
Recursions (\ref{vmf}) and (\ref{vmf0}) provide
\begin{eqnarray*}
\delta_i&=&\delta_{i-1}+\frac{1}{n}a\delta_{i-1}+\mathfrak{q}_n\sigma\xi_i-
\mathfrak{q}_nA\delta_{i-1}+\big(f(t_i)-f(t_{i-1})\big)\mathfrak{q}_n
\nonumber\\
&&\quad -\mathfrak{b}\big(f^{(k)}(t_i)-f^{(k)}(t_{i-1})\big).
\end{eqnarray*}
Multiplying both sides of this equation from the left by $C_n$ we find
\begin{eqnarray}\label{3.3a}
\Delta_i&=&\Delta_{i-1}+\frac{1}{n}C_na\delta_{i-1}
+C_n\mathfrak{q}_n\sigma\xi_i
-C_n\mathfrak{q}_nA\delta_{i-1}
\nonumber\\
&&\quad
+\big(f(t_i)-f(t_{i-1})\big)C_n\mathfrak{q}_n
-C_n\mathfrak{b}\big(f^{(k)}(t_i)-f^{(k)}(t_{i-1})\big).
\end{eqnarray}
A special structure (see \cite{KhL}) of the objects involved in \eqref{3.3a}
\begin{equation*}
\begin{split}
&C_na=n^{\frac{1}{2\beta+1}}aC_n
\\
&C_n\mathfrak{q}_n=n^{-\frac{\beta}{2\beta+1}}\mathfrak{q}
\\
&\frac{1}{n}C_na\delta_{i-1}=n^{-\frac{2\beta}{2\beta+1}}a\Delta_{i-1}
\\
&C_n\mathfrak{q}_nA\delta_{i-1}=n^{-\frac{2\beta}{2\beta+1}}\mathfrak{q}A\Delta_{i-1}
\\
&C_n\mathfrak{b}=n^{\frac{1}{2\beta+1}}\mathfrak{b}
\end{split}
\end{equation*}
allows to simplify \eqref{3.3a} significantly:
\begin{eqnarray*}
\Delta_i&=&\Delta_{i-1}+n^{-\frac{2\beta}{2\beta+1}}\big(a-\mathfrak{q}A\big)\Delta_{i-1}
+n^{-\frac{\beta}{2\beta+1}}\mathfrak{q}\sigma\xi_i
\\
&+&
n^{-\frac{1}{2\beta+1}}\mathfrak{q}\big(f(t_i)-f(t_{i-1})\big)
-n^{\frac{1}{2\beta+1}}\mathfrak{b}\big(f^{(k)}(t_i)-f^{(k)}(t_{i-1})\big).
\end{eqnarray*}
With
$
D_n=I+n^{-\frac{2\beta}{2\beta+1}}\big(a-\mathfrak{q}A)
$
we rewrite the recursion for $\Delta_i$'s to
\begin{eqnarray}\label{3.42}
\Delta_i&=&D_n\Delta_{i-1}
+n^{-\frac{\beta}{2\beta+1}}\mathfrak{q}\sigma\xi_i
\nonumber\\
&+&
n^{-\frac{1}{2\beta+1}}\mathfrak{q}\big(f(t_i)-f(t_{i-1})\big)
-n^{\frac{1}{2\beta+1}}\mathfrak{b}\big(f^{(k)}(t_i)-f^{(k)}(t_{i-1})\big).
\end{eqnarray}
In Proposition 4.1 in \cite{KhL}, it is shown that for $n$ large enough the
magnitudes of all eigenvalues of $D_n$ are strictly less than 1.
We shall use this property for asymptotic analysis and continuous time approximation.

\subsection{Normalized bias and variance}

Denote the bias and variance of the normalized error $\Delta_i$:
$$
M^n_i=E\Delta_i \ \text{and} \ P^n_i=E(\Delta_i-M^n_i)(\Delta_i-M^n_i)^*.
$$
Taking the expectation from both side of
\eqref{3.42} we find
\begin{eqnarray}\label{3.7}
M^n_i&=&D_nM^n_{i-1}
\\
&+&n^{-\frac{1}{2\beta+1}}\mathfrak{q}\big(f(t_i)-f(t_{i-1})\big)
-n^{\frac{1}{2\beta+1}}\mathfrak{b}\big(f^{(k)}(t_i)-f^{(k)}(t_{i-1})\big).
\nonumber
\end{eqnarray}
From \eqref{3.42} and \eqref{3.7}, we get
$
(\Delta_i-M^n_i)=D_n(\Delta_{i-1}-M^n_{i-1})
+n^{-\frac{\beta}{2\beta+1}}\mathfrak{q}\sigma\xi_i,
$
so that $P^n_i=E(\Delta_i-M^n_i)(\Delta_i-M^n_i)^*$ is defined by the
recursion
$$
P^n_i=D_nP^n_{i-1}D_n^*
+n^{-\frac{2\beta}{2\beta+1}}\sigma^2\mathfrak{q}\mathfrak{q}^*.
$$
Since
$
|f^{(k)}(t_i)-f^{(k)}(t_{i-1})|\le \frac{L}{n},
$
it is natural to choose $f$ with
$$
f^{(k)}(t_i)-f^{(k)}(t_{i-1})\equiv\frac{L}{n} \ \text{or} \ -\frac{L}{n}
$$
and substitute
$
\sup_{f\in\sum(\beta,L)}M^n_i\big(M^n_i\big)^*
$
by
$
\widetilde{M}^n_i\big(\widetilde{M}^n_i\big)^*.
$
Henceforth, $\widetilde{M}^n_i$ is defined by recursion (\ref{3.7}) with such $f$, i.e.
$$
\widetilde{M}^n_i=\widetilde{M}^n_{i-1}+n^{-\frac{2\beta}{2\beta+1}}\big(a-
\mathfrak{q}A\big)\widetilde{M}^n_{i-1}+n^{-\frac{2\beta}{2\beta+1}}\varrho^n
-n^{-\frac{2\beta}{2\beta+1}}L\mathfrak{b},
$$
where $\varrho^n=\frac{L\mathfrak{q}}{(k+1)!n^{k}}$.
Thus,
$
\trace\big(P^n_i+\widetilde{M}^n_i\big(\widetilde{M}^n_i\big)^*\big)
$
determines the normalized mean square tracking error
\begin{equation}\label{main01til}
C(\mathfrak{q})=\varlimsup_{n\to\infty}
\trace\Big(P^n_i+\widetilde{M}^n_i\big(\widetilde{M}^n_i\big)^*\Big),
\ t_i\ge c(\mathfrak{q})n^{-\frac{1}{2\beta+1}}\log n.
\end{equation}

\subsection{Continuous time approximation}

To find ``$\varlimsup_{n\to\infty}$'' in \eqref{main01til},
we give
a continuous time approximation of $(\widetilde{M}^n_i,P^n_i)$. To this end, let us
introduce the time stretching
$$
(t_i-t_{i-1})= n^{-1} \ \Rightarrow \ (s_i-s_{i-1})= n^{-\frac{2\beta}{2\beta+1}},
$$
with $t_0=s_0=0$. The boundary layer
$[0,c(\mathfrak{q})n^{-\frac{1}{2\beta+1}}\log n]$ is transformed to
$[0,c(\mathfrak{q})\log n]$ and the interval $[0,T]$ to $[0,Tn^{\frac{1}{2\beta+1}}]$.
Let us define  $\widetilde{M}^n_{s_i}=
\widetilde{M}^n_i$ and $P^n_{s_i}=P^n_i$, $i=0,1,\ldots,$ and for  $s\in [s_{i-1},s_i)$
\begin{eqnarray}\label{3.15'}
\widetilde{M}^n_s&=&\widetilde{M}^n_{s_{i-1}}+\int_{s_{i-1}}^s
\Big((a-\mathfrak{q}A)\widetilde{M}^n_{s_{i-1}}\Big)
ds'+(\varrho^n-\mathfrak{b}L)(s-s_{i-1})
\notag\\
P^n_s&=&P^n_{s_{i-1}}+\int_{s_{i-1}}^s\Big(\big(a-\mathfrak{q}A\big)P^n_{s_{i-1}}+
P^n_{s_{i-1}}\big(a-\mathfrak{q}A\big)^*
\notag\\
&&\hskip .4in+n^{-\frac{\beta}{2\beta+1}}\big(a-\mathfrak{q}A\big)
P^n_{s_{i-1}}\big(a-\mathfrak{q}A
\big)^*+\sigma^2\mathfrak{q}\mathfrak{q}^*\Big)ds'.
\end{eqnarray}
We shall consider these recursions for $s_i$ from
$
\big(c(\mathfrak{q})\log n, \ n^{\frac{1}{2\beta+1}}T\big],
$
where $(\widetilde{M}^n_s,P^n_s)$ have entries bounded in $n$
(see, \cite{KhL}).
Taking into account that recursions (\ref{3.15'}) are homogeneous in $s$,
let us replace $\{s$ and $s_i\}$ by
$\{u=s-c(\mathfrak{q})\log n$ and $u_i=s_i-c(\mathfrak{q})\log n\}$. Then, the entries of
$(\widetilde{M}^n_u,P^n_u)$ are bounded in $n$ for $0\le u\le
n^{\frac{1}{2\beta+1}}T-c(\mathfrak{q})\log n$. Therefore without loss of generality
we may consider (\ref{3.15'}) with initial conditions bounded in $n$:
\begin{equation}\label{3.17'}
C(\mathfrak{q})=\lim_{u\to\infty}\varlimsup_{n\to\infty}
\trace\Big(P^n_u+\widetilde{M}^n_u\big(\widetilde{M}^n_u\big)^*\Big).
\end{equation}
To determine the right hand side of (\ref{3.17'}), we apply the Arzela-Ascoli theorem.
For any $T>0$, the functions $(\widetilde{M}^n_u,P^n_u)_{0\le u\le T}$ are uniformly
bounded and equicontinuous.
So, by the Arzela-Ascoli theorem, any converging subsequence
$(\widetilde{M}^{n'}_u,P^{n'}_u)$
obeys the limit $(\widetilde{M}'_u,P'_u)$ in the local uniform topology:
\[
\lim_{n'\to\infty}\sum_{N=1}^\infty\frac{1}{2^N}
\min\Big(1,\sup_{0\le u\le N}\big(\|\widetilde{M}^{n'}_u-\widetilde{M}'_u\|+
\|P^{n'}_u-P'_u\|\big)\Big)=0,
\]
where
\begin{eqnarray*}
&&
\widetilde{M}'_u=\widetilde{M}'_0+\int_0^u\Big(\big(a-\mathfrak{q}A\big)
\widetilde{M}'_v+\mathfrak{b}L\Big)dv
\\
&&P'_u=P'_0+\int_0^u\Big(\big(a-\mathfrak{q}A\big)P'_v+P'_v\big(a-\mathfrak{q}A\big)^*+
\sigma^2\mathfrak{q}\mathfrak{q}^*\Big)dv.
\end{eqnarray*}
Since the eigenvalues of $a-\mathfrak{q}A$ have negative real parts, the limits
$\widetilde{M}(\mathfrak{q}):=\lim_{u\to\infty}\widetilde{M}^{'}_u$
and
$P(\mathfrak{q}):=\lim_{u\to\infty}P'_u$ exist and are defined as:
\begin{eqnarray}\label{3.8}
&&\widetilde{M}(\mathfrak{q})=-
L\big(a-\mathfrak{q}A\big)^{-1}\mathfrak{b}
\\
&&
\big(a-\mathfrak{q}A\big)P(\mathfrak{q})+P(\mathfrak{q})\big(a-\mathfrak{q}A\big)^*+
\sigma^2\mathfrak{q}\mathfrak{q}^*=0,
\end{eqnarray}
that is $\widetilde{M}(\mathfrak{q})$ and
$P(\mathfrak{q})$ are independent of $\{n'\}$ and so
\begin{equation}\label{prosto}
C(\mathfrak{q})=\trace\big(P(\mathfrak{q})+\widetilde{M}(\mathfrak{q})
\widetilde{M}^*(\mathfrak{q})\big).
\end{equation}

\section{Minimization of the cost function in one parameter class}
\label{sec-4}

\subsection{Motivation}

For large values of $k$ a direct minimization of
$
C(\mathfrak{q})
$
from \eqref{prosto} would be a difficult problem. Moreover,
$$
\mathfrak{q}^\circ=\argmin_\mathfrak{q}C(\mathfrak{q})
$$
could not a priori guarantee  negative real parts of eigenvalues for
$(a-\mathfrak{q}^\circ A)$. To avoid implementation of a conditional minimization
procedure, we propose to choose $\mathfrak{q}$ from some
limited class given below.

\subsection{Adaptation to Kalman filter design}

Our estimator has a structure of Kalman filter in the discrete time.
We assume that $F(0)$ is a random vector, $ \widehat{F}^n(0)=
EF(0)$, and $f^{(k)}(t_i)$ is generated by stochastic
recursion
\begin{equation}\label{gl1}
f^{(k)}(t_i)=f^{(k)}(t_{i-1})+ n^{-\frac{\beta+1}{2\beta+1}}\gamma\eta_i,
\end{equation}
where $(\eta_i)$ is a white noise, independent of $(\xi_i)$, with $E\eta_1=0$,
$E\eta_1=1$ and $\gamma$ is an arbitrary nonzero parameter. For the observation model
$$
X_i=f^{(0)}(t_{i-1})+\sigma\xi_i
$$
we apply the estimator given in (\ref{vmf}). The resulting errors
$\delta_i=\widehat{F}(t_i)-F(t_i)$, $i=1,\ldots,$ are defined by a recursion
$$
\delta_i=\delta_{i-1}+\frac{1}{n}a\delta_{i-1}+\mathfrak{q}_n\sigma\xi_i-
\mathfrak{q}_nA\delta_{i-1}-n^{-\frac{\beta+1}{2\beta+1}}\mathfrak{b}\gamma \eta_i.
$$
Then, for $\Delta_i=C_n\delta_i, \ i\ge 1$
we obtain
\begin{equation}\label{3.411}
\Delta_i=\Delta_{i-1}+n^{-\frac{2\beta}{2\beta+1}}\big(a-\mathfrak{q}A\big)\Delta_{i-1}
+n^{-\frac{\beta}{2\beta+1}}\mathfrak{q}\sigma\xi_i
-n^{-\frac{\beta}{2\beta+1}}\mathfrak{b}\gamma \eta_i
\end{equation}
and supply $\Delta_0=\delta_0$.
Denote $Q^n_i=E\Delta_i\Delta^*_i$. From (\ref{3.411}) it follows
\begin{eqnarray*}
Q^n_i&=&\Big(I+n^{-\frac{2\beta}{2\beta+1}}\big(a-\mathfrak{q}A\big)\Big)Q^n_{i-1}
\Big(I+n^{-\frac{2\beta}{2\beta+1}}\big(a-\mathfrak{q}A\big)\Big)^*
\\
&&+n^{-\frac{2\beta}{2\beta+1}}\sigma^2\mathfrak{q}\mathfrak{q}^*
+n^{-\frac{2\beta}{2\beta+1}}\gamma^2\mathfrak{b}\mathfrak{b}^*
\\
&=&D_nQ^n_{i-1}D^*_n+n^{-\frac{2\beta}{2\beta+1}}\sigma^2\mathfrak{q}\mathfrak{q}^*
+n^{-\frac{2\beta}{2\beta+1}}\gamma^2\mathfrak{b}\mathfrak{b}^*.
\end{eqnarray*}
Similar to $(\widetilde{M}^n_s, \ P^n_s)$ (see previous section), let us introduce
$Q^n_s$:
\begin{multline*}
Q^n_s=Q^n_{s_{i-1}}+\int_{s_{i-1}}^s\Big(\big(a-\mathfrak{q}A\big)Q^n_{s_{i-1}}+
Q^n_{s_{i-1}}\big(a-\mathfrak{q}A\big)^*
\\
+n^{-\frac{\beta}{2\beta+1}}\big(a-\mathfrak{q}A\big)Q^n_{s_{i-1}}\big(a-\mathfrak{q}A
\big)^*+\sigma^2\mathfrak{q}\mathfrak{q}^*+\gamma^2\mathfrak{b}\mathfrak{b}^*\Big)ds'.
\end{multline*}
Applying  the Arzela-Ascoli theorem technique it can be readily shown that
$Q^{n}_s$ converges in the local uniform topology to $Q_s$,
where
$$
Q_s=Q_0+\int_0^s\Big(\big(a-\mathfrak{q}A\big)Q_{s'}+Q_{s'}\big(a-\mathfrak{q}A\big)^*+
\sigma^2\mathfrak{q}\mathfrak{q}^*+\gamma^2\mathfrak{b}\mathfrak{b}^*\Big)ds,
$$
and $\lim_{s\to\infty}Q_s:=Q$ with $Q$ being the unique solution of Lyapunov
equation
\begin{equation}\label{Lyap}
\big(a-\mathfrak{q}A\big)Q+Q\big(a-\mathfrak{q}A\big)^*+\sigma^2\mathfrak{q}
\mathfrak{q}^*+\gamma^2\mathfrak{b}\mathfrak{b}^*=0.
\end{equation}
The matrix $Q$ is a function of arguments $\mathfrak{q}$ and $\gamma$:
$Q=Q(\mathfrak{q},\gamma)$. We choose $\mathfrak{q}=\mathfrak{q}(\gamma)$ so that
for any $\gamma$
\begin{equation}\label{kkk}
Q(\mathfrak{q},\gamma)\ge Q(\mathfrak{q}(\gamma),\gamma):=Q(\gamma)>0.
\end{equation}
Due to the Kalman filtering theory, the lower bound \eqref{kkk} holds true
for
\begin{equation}\label{qg}
\mathfrak{q}(\gamma)=\frac{Q(\gamma)A^*}{\sigma^2}
\end{equation}
with $Q(\gamma)$ being solution of the algebraic Riccati equation
\begin{equation}\label{Riccati}
aQ(\gamma)+Q(\gamma)a^*+\gamma^2\mathfrak{b}\mathfrak{b}^*-\frac{Q(\gamma)A^*AQ(\gamma)}
{\sigma^2}=0.
\end{equation}
It is well known (see e.g. Theorem 16.2 in \cite{LSII}) that \eqref{Riccati} possesses
a unique positive-definite solution
provided that block-matrices
$$
G_1=  \begin{pmatrix}
    A \\
    Aa \\
    \vdots\\
    Aa^k
  \end{pmatrix}
  \quad\text{and}\quad
 G_2= \begin{pmatrix}
    \mathfrak{b}\mathfrak{b}^* & a\mathfrak{b}\mathfrak{b}^* & \ldots &
    a^k\mathfrak{b}\mathfrak{b}^*
  \end{pmatrix}
$$
have full ranks $r=k+1$. Notice that $G_1$ is a unite matrix and the rank of $G_2$ is
$k+1$. Consequently, the eigenvalues of the matrix $(a-\mathfrak{q}(\gamma)A)$
with $\mathfrak{q}(\gamma)$ defined in (\ref{qg}) have negative real parts
(see, Lemma 16.11 in \cite{LSII}).

\subsection{Minimization of the cost function}
\label{sec-4.3}
The one parameter  family
$$
\Gamma=\Big\{\gamma\ge \gamma_\varepsilon>0:\mathfrak{q}(\gamma)=\frac{Q(\gamma)A^*}
{\sigma^2}\Big\}
$$
permits a simple numeric implementation and guarantees filtering stability
mentioned above. In this class we use a constrain parameter
$\gamma_\varepsilon$ to  compensate unacceptably large boundary layer when
the minimization procedure yields small values of
$
\gamma^\circ=\argmin_{\gamma>0}C\big(\mathfrak{q}(\gamma)\big),
$
with $C\big(\mathfrak{q}(\gamma)\big)$ given in \eqref{prosto}.
The magnitude of $\gamma_\varepsilon$ is dictated by $k$, $\sigma$,
initial conditions and boundary layer specifications.
The minimization in our class provides
\begin{equation}\label{choice}
\gamma^\circ=\mathop{{\rm argmin}}_{\gamma\ge \gamma_\varepsilon}
C\big(\mathfrak{q}(\gamma)\big)\quad\text{and}\quad
\mathfrak{q}(\gamma^\circ)=\frac{Q(\gamma^\circ)A^*}{\sigma^2}.
\end{equation}

\section{Explicit minimization procedure}
\label{sec-5}

\subsection{Filter gain}
In this section we describe a structure of $\mathfrak{q}(\gamma^\circ)$.
Recall that $\mathfrak{q}(\gamma)=\frac{Q(\gamma)A^*}{\sigma^2}$ and $Q(\gamma)$ solves
the Riccati equation \eqref{Riccati} for any fixed $\sigma$. For notational
convenience replace $Q(\gamma)$ by $Q(\gamma,\sigma)$ and set
$
U=Q(1,1).
$
Clearly, $U$ solves the algebraic Riccati equation
\[
aU+Ua^*+\mathfrak{b}\mathfrak{b}^*-UA^*AU=0.
\]
Kalachev \cite{K} shows that
$$
Q_{ij}(\gamma,\sigma)=U_{ij}\sigma^2
\Big({\frac{\gamma}{\sigma}}\Big)^{\frac{i+j+1}{k+1}}, \ i,j=0,1,\ldots,k,
$$
where $Q_{ij}(\gamma,\sigma)$ and $U_{ij}$ are entries of $Q(\gamma,\sigma)$ and
$U$ respectively.
Hence,
\begin{eqnarray*}
&&
q_0(\gamma)=U_{00}\Big(\frac{\gamma}{\sigma}\Big)^{1/k+1}
\\
&&
q_1(\gamma)=U_{01}\Big(\frac{\gamma}{\sigma}\Big)^{2/k+1}
\\
&&.................................
\\
&&
q_k(\gamma)=U_{0k}\Big(\frac{\gamma}{\sigma}\Big).
\end{eqnarray*}
For $k\le 4$, these values are given in the table below.

\medskip
$$
\mbox{
\begin{tabular}{|c|c|c|c|c|c|}
\hline
  k & $U_{00}$ & $U_{01}$ & $U_{02}$ & $U_{03}$ & $U_{04}$ \\
\hline
  0 & 1 & NA & NA & NA & NA \\
  1 & $\sqrt{2}$ & 1 & NA & NA & NA \\
  2 & 2 & 2 & 1 & NA & NA \\
  3 & $\sqrt{4+\sqrt{8}}$ & $2+\sqrt{2}$ & $\sqrt{4+\sqrt{8}}$ & 1 & NA \\
  4 & $1+\sqrt{5}$ & $3+\sqrt{5}$ & $3+\sqrt{5}$ & $1+\sqrt{5}$ & 1 \\ \hline
\end{tabular}
}
\eqno(\rm{Table} \ 1)
$$

\begin{figure}
\begin{center}
\includegraphics[scale=0.7]{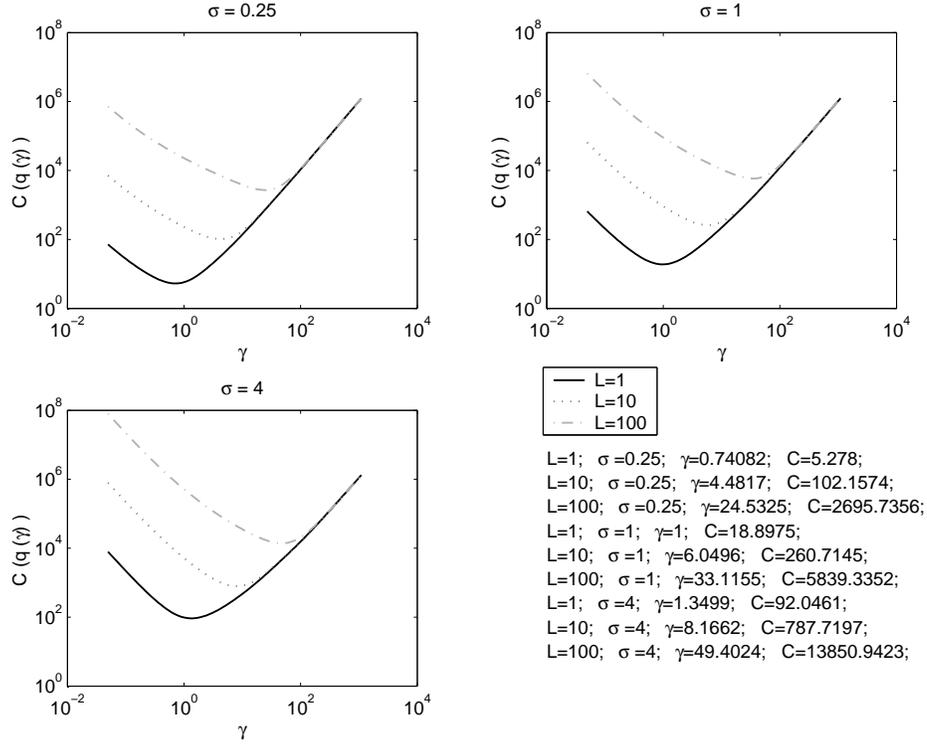}
\end{center}
\caption{Cost function $C\big(\mathfrak{q}(\gamma)\big)$ in logarithmic scale
for various $L$ and $\sigma$; $k=2$.
}
\label{fig1}
\end{figure}

\medskip
The complex structure of $C\big(\mathfrak{q}(\gamma)\big)$ does
not provide an insight of  $\gamma$ and $L$ connection. Numerical simulations
show that for a wide
range of values  $\log(\gamma^\circ)$ is almost
proportional to $\log(L)$ (see also Figure \ref{fig1}).  This remark enables to construct a simple
interpolation tables for the values of $\gamma^\circ$ and
$C\big(\mathfrak{q}(\gamma))$ with respect to the parameter $L$.

\subsection{Eigenvalues of $\b{(a-\mathfrak{q}(\gamma)A)}$}

Although the  eigenvalues of $\mathfrak{q}(\gamma)$ have negative real parts,
we may not formally guarantee that they are different.
So, for $k\le 4$ we give the eigenvalues:
\begin{eqnarray*}
&&
k=0: \quad-\Big(\frac{\gamma}{\sigma}\Big)
\\
&&
k=1: \quad-\Big(\frac{\gamma}{\sigma}\Big)^{1/2}
\Big(\frac{1}{\sqrt{2}}\pm i\frac{1}{\sqrt{2}}\Big)
\\
&&
k=2: \quad-\Big(\frac{\gamma}{\sigma}\Big)^{1/3}\Big(1; \ \frac{1}{2}
\pm i\frac{\sqrt{3}}{2}\Big)
\\
&&
k=3: \quad-\Big(\frac{\gamma}{\sigma}\Big)^{1/4}\Big(0.924\pm i
0.383; \ 0.383\pm i0.924\Big)
\\
&&
k=4:
\quad-\Big(\frac{\gamma}{\sigma}\Big)^{1/5}\Big(1; \ 0.809\pm i
0.588;  \ 0.309\pm i0.951\Big).
\end{eqnarray*}
For $k>4$, the fact that the polynomial has different roots should be verified.

Notice that $(\gamma/\sigma)^{1/(k+1)}$ is a natural
control parameter  defining the size of
boundary layer and should be limited from below by
$(\gamma_\varepsilon/\sigma)^{1/(k+1)}$ with appropriate $\gamma_\varepsilon$.

\subsection{Cost function}
The index form of
$\big(a-\mathfrak{q}A\big)\widetilde{M}(\mathfrak{q})=-L\mathfrak{b}$
i.e.
$$
  \begin{pmatrix}
-q_0 & 1 & 0 & 0 &\cdots & 0
\\
-q_1 & 0 & 1 & 0 &\cdots & 0
\\
\vdots & \vdots & \vdots & \vdots &\cdots & \vdots
\\
-q_{k-1} & 0 & 0 & 0 &\cdots & 1
\\
-q_k & 0 & 0 & 0 &\cdots & 0
  \end{pmatrix}
  \begin{pmatrix}
    \widetilde{M}_0(\mathfrak{q})
    \\
    \widetilde{M}_1(\mathfrak{q})
    \\
    \vdots
    \\
    \widetilde{M}_{k-1}(\mathfrak{q})
    \\
    \widetilde{M}_k(\mathfrak{q})
  \end{pmatrix}=
  -L
  \begin{pmatrix}
    0
    \\
    0
    \\
    \vdots
    \\
    0
    \\
    1
  \end{pmatrix}
$$
allows to find the solution
$$
\widetilde{M}_0(\mathfrak{q})=L\frac{1}{q_k},\widetilde{M}_1(\mathfrak{q})=
L\frac{q_0}{q_k},\ldots,\widetilde{M}_k(\mathfrak{q})=L\frac{q_{k-1}}{q_k}.
$$
Thus,
$
\trace \Big(\widetilde{M}(\mathfrak{q})\widetilde{M}^*(\mathfrak{q})\Big)
=L^2
\Big[\Big(\frac{1}{q_k}\Big)^2+\sum_{j=0}^{k-1}
\Big(\frac{q_j}{q_k}\Big)^2\Big].
$
From \eqref{3.8} it follows $
P(\mathfrak{q})=\sigma^2\int_0^\infty\mathfrak{q}^*e^{(a-\mathfrak{q}A)^*t}
e^{(a-\mathfrak{q}A)t}\mathfrak{q}dt$.

As a result, the final expression for the cost function is
$$
C\big(\mathfrak{q}\big)=\sigma^2\Bigg(\int_0^\infty\mathfrak{q}^*e^{(a-\mathfrak{q}A)^*t}
e^{(a-\mathfrak{q}A)t}\mathfrak{q}dt+\Big(\frac{L}{\sigma}\Big)^2
\Big[\Big(\frac{1}{q_k}\Big)^2+\sum_{j=0}^{k-1}
\Big(\frac{q_j}{q_k}\Big)^2\Big]\Bigg).
$$

\section{Conclusion remark}

In this paper, we use the fact that a class of Kalman filters, being adapted  to a
nonparametric statistic setting, provides the optimal rate of convergence in sample
size ($n\to\infty$).
We show how to evaluate a normalized risk function for large sample size
and minimize that value in some subclass of Kalman filters with constant filter gain.
The Kalman type estimator, as any on-line estimator, has inevitable boundary layer.
We suggest to reduce the boundary layer by interpolation procedure and limitation
from below for filtering gain.


\vskip .4in
\begin{thebibliography}{9}

\bibitem{B_Sh} Bar-Shalom, Yaakov and Li, Xiao-Rong {\em Estimation and tracking.
Principles, techniques, and software.} Artech House, Inc., Boston, MA, 1993.

\bibitem {CKL}  Chow, P.-L.,  Khasminskii, R. and  Liptser, R.Sh. (1997)
Tracking of signal and its derivatives in Gaussian white noise, {\it Stochastic
processes and their application.} {\bf 69}, 2  pp. 259-273.

\bibitem{IK80} Ibragimov, I. and Khasminskii, R.  (1980) On nonparametric
estimation of regression,{\em Soviet Math.Dokl.},{\bf 21}, pp. 810--814.

\bibitem{IK81}  Ibragimov, I.  and  Khasminskii, R.  {\em Statistical estimation:
Asymptotic theory.} Springer Verlag, 1981 (Russian ed.1979).

\bibitem{K}  M.G. Kalachev ``One method of multiple differentiation applied
to a signal in automatic regulation systems''. {\em Automation and Remote
Control}, No 6, 1970, pp. 890--896.

\bibitem{KhL} Khasminskii, R and Liptser, R. (2001) On-line estimation of a smooth
regression function

\bibitem{LSII}  R.Sh. Liptser and A.N. Shiryayev {\em Statistics of Random
Processes, II.} Springer Verlag. 2000

\bibitem{P}  Parzen, E. (1962) On estimation of a probability density function and
mode. {\em Ann. Math. Statist.}, {\bf 33}, No 3, pp. 1065--1073.

\bibitem{R}  Rosenblatt, M. (1956) Remarks on some nonparametric estimates of a
density function''. {\em Ann. Math. Statist.}, {\bf 27}, No 3, pp. 832--837.

\bibitem{St80} Stone, C. (1980) Optimal rates of convergence for
nonparametric estimators, {\em Ann. Statist.}, {\bf 8}, pp.1348--1360.

\bibitem{St82} Stone, C.  (1982) Optimal global rates of convergence for
nonparametric regression, {\em Ann. Statist.}, {\bf 10}, pp.1040--1053.
\end{thebibliography}
\end{document}